\documentclass[12pt,reqno]{amsart}

\usepackage{graphicx,subfigure}
\usepackage{hyperref}
\hypersetup{colorlinks=true, citecolor=blue, linkcolor=red}

\numberwithin{equation}{section} \numberwithin{figure}{section}
\numberwithin{table}{section} 
{ \theoremstyle{plain}

}

{ \theoremstyle{definition}

\numberwithin{equation}{section} \numberwithin{lem}{section}
\numberwithin{thm}{section} \numberwithin{cor}{section}
\numberwithin{pro}{section} \numberwithin{rem}{section}
}

\begin{document}

\title{A small remark on the derivation of the Plateau angle conditions for the vector-valued Allen-Cahn equation}

\author{Christos Sourdis}
\address{Department of Applied Mathematics and Department of Mathematics, University of
Crete, 700 13 Panepistimioupoli Vouton, Crete, Greece.}
\email{csourdis@tem.uoc.gr} \maketitle
\begin{abstract}
We clarify a point in \cite{damialis}, concerning the derivation of
the Plateau angle conditions for the vector-valued Allen-Cahn
equation.
\end{abstract}
The Plateau angle conditions for the vector-valued Allen-Cahn
equation have recently been  rigorously derived in the interesting
preprint \cite{damialis}. However, it is not clear to the author how
the estimates (30) and (31) in \cite{damialis} are used for showing
that the terms involving integrals of $T_{ij}(v)\frac{y_j}{y_2}$
with $j\neq 2$  tend to zero as $R\to \infty$ in (27) of the latter
reference.
 In this note we will clarify this point and verify that the aforementioned terms indeed tend to
 zero.
 To this end, instead of bounding the absolute value of these
integrals (as was done in \cite{damialis}), we will exploit a
cancelation property that seems to have been left unnoticed in
\cite{damialis}.

Throughout this short note we will follow exactly the notation of
\cite{damialis}.

In what follows, we will show that
\[
\frac{1}{R}\int_{-R\sin \psi_1}^{R\sin \psi_1}\int_{-R\cos
\psi_2}^{R\cos \psi_2} T_{33}(v)\frac{y_3}{y_2}dy_3 dy_1 \to 0\ \
\textrm{as}\ R\to \infty
\]
(the remaining terms can be handled analogously).

We note that
\[
T_{33}(v)=\frac{1}{2}\left(|v_{,3}|^2-|v_{,1}|^2-|v_{,2}|^2-2W(v)\right)\
\ (v_{,i}=\partial v/ \partial y_{i}),
\]
and $y_2=\sqrt{R^2-y_1^2-y_3^2}$. Moreover, we point out that
$R\sin\psi_{1}(R)\to \infty$ and $\psi_{2}(R)\to 0$ as $R\to
\infty$.
 By virtue of Lebesgue's dominated
convergence theorem, whose assumptions have been verified in
\cite{damialis}, it suffices to show that, given $y_1 \in
\mathbb{R}$, we have
\[
\frac{1}{R}\int_{-R\cos \psi_2}^{R\cos \psi_2}
T_{33}(v)\frac{y_3}{y_2}dy_3  \to 0\ \ \textrm{as}\ R\to \infty.
\]
Letting $\tilde{y}_3=R^{-1}y_3$, with a slight abuse of notation, we
can write the above integral as
\[
\frac{1}{2}\int_{-\cos \psi_2}^{\cos \psi_2}
\left(\frac{1}{R^2}|v_{,3}|^2-|v_{,1}|^2-|v_{,2}|^2-2W(v)
\right)\frac{R \tilde{y}_3}{\sqrt{R^2-y_1^2-R^2
\tilde{y}_3^2}}d\tilde{y}_3.
\]
Now, by Hypothesis 2 in \cite{damialis}, and a trivial application
of Lebesgue's dominated convergence theorem, we infer that, as $R\to
\infty$, the above integral converges to
\[
-\left(\frac{1}{2}\left|\dot{U}_{12}(y_1) \right|^2+W\left(
U_{12}(y_1)\right)\right)\int_{-1}^{1} \frac{ \tilde{y}_3}{\sqrt{1-
\tilde{y}_3^2}}d\tilde{y}_3=0,
\]
as desired.

\end{document}